\documentclass[twocolumn]{autart}    
\usepackage{amsmath,amssymb}
\usepackage{color}

\usepackage{xargs, xcolor}
\newcommand{\Sn}{\mathbb{S}_+^n}
\newcommand{\coneK}{\mathcal{K}}
\newcommand{\vsV}{\mathcal{V}}
\newcommand{\excdZero}{\backslash \{0 \}}
\newcommand{\dH}{d_{\mathcal{H}}}

\usepackage{wrapfig}

\newtheorem{definition}{Definition}

\newtheorem{remark}{Remark}

\newtheorem{lemma}{Lemma}
\newtheorem{fct}{Fact}
\newtheorem{corollary}{Corollary}
\usepackage{dsfont}
\newcommand{\allOne}{\mathds{1}}

\usepackage{tcolorbox}

\begin{document}

\begin{frontmatter}

\title{Lyapunov and Riccati Equations from a Positive System Perspective} 

\thanks[footnoteinfo]{Corresponding author: Yankai~Lin. 
This project has received funding from the European Research Council
(ERC) under the European Union’s Horizon 2020 research and innovation
programme under grant agreement No 834142 (ScalableControl).
}

\author[DJW]{Dongjun Wu}\ead{dongjun.wu@control.lth.se},    
\author[YKL]{Yankai Lin}\ead{25021103@wit.edu.cn} 

\address[DJW]{Department of Automatic Control, Lund University, Sweden}  
\address[YKL]{School of Computer Science and Engineering, Wuhan Institute of Technology, Wuhan, PR China}

\begin{abstract}                          
This paper presents a new interpretation of the Lyapunov and Riccati equations from the perspective of positive system theory.
We show it is possible to construct positive systems related to these equations, and then certain conclusions --  such as the existence and uniqueness of solutions -- can be drawn from positive systems theory. Specifically, under standard observability assumptions, a strictly positive linear system can be constructed for Lyapunov equations, leading to exponential convergence in Hilbert metric to the Perron-Frobenius vector -- closely related to the solution of the Lyapunov equation.
For algebraic Riccati equations, homogeneous strictly positive systems can be constructed, which exhibit more complex dynamical behaviors. While the existence and uniqueness of the solution can still be proven, only asymptotic convergence can be obtained.

\end{abstract}

\end{frontmatter}

\section{Introduction }

The Lyapunov and algebraic Riccati equations play fundamental roles in control
theory, underpinning applications in system analysis, model order reduction, stabilization
and optimal control \cite{datta2004numerical,ito2025weighted,moore1981principal,abou2012matrix,lancaster1995algebraic}. 
Their properties -- such as existence, uniqueness, and numerical
computation of solutions -- have been extensively studied in the
literature. 
In particular, it is well-known that Lyapunov equations for linear time-invariant systems admit explicit solutions - expressed as integrals for continuous-time systems or infinite sums for discrete-time systems. In contrast, algebraic Riccati equations rarely yield such explicit solutions, leading to a distinct treatment of these equations despite their strong connections.

Recently, Pates \cite{pates2024lyapunov} proposed viewing the Lyapunov equation as a fixed-point problem on a compact convex set, leveraging the Brouwer fixed-point theorem to establish solution existence and prove uniqueness under an observability assumption. However, this technique does not directly extend to Riccati equations. 
To address this, we propose a new approach to tackle Lyapunov
and algebraic Riccati equations in a unified framework, namely, by associating these equations with dynamical systems, and the properties of the dynamical systems in turn characterize the solutions of these equations.

Our main research finding is that both the Lyapunov and algebraic Riccati
equations can be analyzed through the lens of positive systems theory
-- a framework that has emerged as a powerful tool in control theory
\cite{rantzer2018tutorial,lam2019positive}.
To be specific, we first consider the simpler case of Lyapunov equations in both continuous-time and discrete-time. By introducing an auxiliary variable, we construct a linear matrix differential (difference) equation which is shown to be strictly positive under standard observability conditions. Following standard arguments from positive systems theory, we demonstrate exponential convergence of the system in Hilbert metric to the unique Perron-Frobenius vector, which further leads to the unique solution to the Lyapunov equation. We then apply similar techniques to algebraic Riccati equations, leading to homogeneous strictly positive systems. 
In such cases, we prove asymptotic convergence (in Hilbert metric) of a dynamical system to the unique solution of the algebraic Riccati equation.

Our contributions are twofold:
\begin{enumerate}
    \item We establish novel connections between positive systems theory and classical equations in control, namely Lyapunov and algebraic Riccati equations. A key technical step involves homogenizing these equations by introducing a scalar-valued linear function, which enables a unified treatment under the framework of positive systems.

    \item Leveraging both linear and nonlinear positive systems theory, we characterize the solutions of Lyapunov and algebraic Riccati equations as steady states of suitably constructed positive dynamical systems. Furthermore, we analyze the convergence properties of these dynamics, thereby providing a dynamical systems perspective on solving such equations.
\end{enumerate}

The paper is organized as follows. In Section \ref{sec:preliminary}, we recall some concepts from positive systems theory and state some useful results that will be used as the main tools in this paper. In Section \ref{sec:lyap}, we study Lyapunov equations by positive systems theory. Section \ref{sec:ricc} extends the results in Section \ref{sec:lyap} to algebraic Riccati equations.

{\it Notations:} Given a set $S$, the set of interior points of $S$ is denoted as ${\rm Int} \, S$. For a square matrix $A$, the trace of $A$ is denoted as ${\rm tr}(A)$. 
Let $A_1, \cdots , A_m$ be some matrices, denote
${\rm col}(A_1, \cdots ,A_m)$ as the column stack of these matrices.
Given two matrices $A\in \mathbb{R}^{n\times n}$ and $C\in \mathbb{R}^{m\times n}$, we say that $(A,C)$ is observable if ${\rm col}(C, CA,\cdots ,CA^{n-1} ) $ is full rank.
The set of $n\times n$ positive semi-definite matrices is denoted as $\Sn$, and we write $X \succ Y$ ($ X \succeq Y$) if $X-Y$ is positive definite (semi-definite). 
For a square matrix $A$, we say that $A$ is \emph{stable} if the eigenvalues of $A$ lie on the left open-half complex plane;
we say that $A$ is \emph{Schur stable} if the eigenvalues of $A$ lie in the
open unit disk in the complex plane. For an ordinary differential equation
$\dot{x}= f(x)$, the solution at time $t$ with initial condition $(t_0=0, x(0)=x_0)$ is 
denoted $\phi(t,x_0)$.

\section{Preliminaries on Positive Systems Theory} \label{sec:preliminary}
In this section, we shall introduce the necessary notations and recall basic concepts that will be used throughout the paper. We begin by defining the cones, and then, on this basis, briefly review relevant results from the theory of positive systems and fixed-point theorems.
\begin{definition}[Cone]
Let $\mathcal{V}$ be a finite-dimensional vector space equipped with a norm $\| \cdot \|$. The subset $\mathcal{K} \subseteq \mathcal{V}$ is called a \emph{cone} if the following three properties are simultaneously satisfied:
\begin{enumerate}
    \item[(a)] ${\rm Int}\,\mathcal{K}\neq\emptyset$; 
    \item[(b)] $\lambda\in\mathbb{R}_{\ge0}$ and $x,\;y\in \mathcal{K}$ $\Rightarrow$
$\lambda x\in \mathcal{K}$, and $x+y\in \mathcal{K}$;
    \item[(c)] $x\in \mathcal{K}\backslash\{0\}$ $\Rightarrow$ $-x\notin \mathcal{K}$. 
\end{enumerate}

In particular, the cone $\mathcal{K}$ is called a \emph{normal cone} if the following condition is further satisfied:
\begin{enumerate}
    \item[(d)] there exists some $\delta>0$ such that $\|x+y\|>\delta$ whenever
$x,y\in \mathcal{K}$, with $\|x\|=\|y\| = 1$. 
\end{enumerate}
A partial ordering 
can be defined on $\mathcal{K}$: if $x - y \in \mathcal{K}$, we say $ x \leq_{\coneK} y$.
We also write $x\leq y$ when clear from the context.

\end{definition}

Throughout this paper, we shall encounter two types of cones, namely, the
positive orthant $\mathbb{R}_{+}^{n}:=\{x\in\mathbb{R}^{n}:x_{i}\ge0,\;\forall i\in\{1,\cdots,n\}\}$
and the set of semi-positive definite matrices $\mathbb{S}_{+}^{n}=\{X\in\mathbb{R}^{n\times n}:X \succeq 0\}$. The norm on $\Sn$ is defined as the Frobenius norm $\|P \|^2 = \sum_{ij} P_{ij}^2$, and the partial ordering is defined as $X \preceq Y \Leftrightarrow X - Y \preceq 0$.
Obviously, both $\mathbb{R}_{+}^{n}$ and $\mathbb{S}_{+}^{n}$ are normal cones

\begin{definition}[Hilbert metric]
Given a cone $\mathcal{K}$ and $x,y\in \coneK \backslash\{0\}$, let $M(x/y):=\inf\{\lambda:x\le\lambda y\}$ and $m(x/y):=\sup\{\mu:\mu y\le x\}$. Then, the Hilbert metric on $\mathcal{K}$ is defined as 
$$
d_{\mathcal{H}}(x,y)=\ln\frac{M(x/y)}{m(x/y)}.
$$
\end{definition}
The Hilbert metric $d_{\mathcal{H}}(x,y) $ measures the distance between the rays $\{ rx; r \ge 0 \}$
and $\{ ry : r \ge 0 \}$, in particular, $d_{\mathcal{H}}(x,y)=0$ if
and only if $x$ and $y$ lie on the same ray. Hilbert metric plays a crucial
role in positive system theory, for instance in consensus \cite{Wu2025,Sepulchre2010}.
For more details about Hilbert metric, see \cite{Bushell1973}.

\begin{definition}[Order-preserving maps]
An operator $T:\vsV\to \vsV$ is called \emph{homogeneous} if $T(\alpha x)=\alpha T(x)$
for all $x\in V,$ $\alpha\in\mathbb{R}$. 
Let $\coneK\subseteq \vsV$ be a cone. Assume $T\coneK \subseteq \coneK$.
We say that $T$ is \emph{order-preserving} if 
\[
0 \le_{\coneK} x \le_{\coneK} y  \; \Rightarrow \;  0 \le_\coneK T  x \le_{\coneK} Ty
\]
We say $T$ is  \emph{strictly order-preserving} if 
\[
0\le_{\coneK} x \le_{\coneK} y \text{ and $x\ne y$ } 
\; \Rightarrow \;
Ty - Tx \in {\rm Int} \; \coneK.
\]
Let $T$ be a homogeneous, continuous and order-preserving
mapping, the \emph{cone spectral radius} of $T$ is defined as
\[
\rho(T)=\sup_{x\in \coneK}\limsup_{n\to\infty}\|T^{n}x\|^{1/n}.
\]
\end{definition}

The following theorem collects some important facts from nonlinear Perron-Frobenius theory.
\begin{thm}
\label{thm:main-pos}Let $ \coneK$ be a normal cone in a finite dimensional vector
space $\vsV$. 
Let $F: \vsV \to \vsV $ be a continuous, homogeneous and order-preserving map on $\coneK$.
Then there exists an $x_*\in \coneK \backslash \{0 \}$, such that $F(x_{*}) = \rho(F) x_*$.
In addition,
\begin{enumerate}
    \item[(a)] if $F$ is strongly order-preserving, then 
    $x_* \in {\rm Int} \; \coneK$ and it is unique (up to positive scaling). Moreover, 
    \begin{equation*}
        d_{\mathcal{H}}(F^k(x),x_{*}) \to 0\text{ as }k\to\infty\text{ for all }x 
        \in \coneK \backslash \{0 \};
    \end{equation*}
    \item[(b)] if $F$ is linear and strongly order-preserving, then there exists a positive constant $\lambda \in [0,1)$, such that 
    \begin{equation*}
        d_{\mathcal{H}}(F^k(x),x_{*}) \le \lambda^k d_{\mathcal{H}}(x,x_{*}), \; 
    \end{equation*}for all $ k\ge 1, \, x \in \coneK \backslash \{0\}$.
\end{enumerate}

\end{thm}

\begin{remark}
    For the existence of $x_*$, see \cite[Corollary 5.4.2]{lemmens2012nonlinear}. 
    For (a), see
    \cite[Theorem 6.5.1]{lemmens2012nonlinear}. In particular, it implies the uniqueness of $x_*$.
    For (b), we refer to \cite{Bushell1973}.
\end{remark}

With these preparations in place, we now define the central objects of the note, namely, positive systems.
\begin{definition}[Positive systems]
\label{def:pos-sys}
Let $ \coneK$ be a normal cone in a finite dimensional vector
space $\vsV$, and $F$ be a continuous, homogeneous mapping on $\vsV$.
Consider the dynamics
\begin{equation}
   \label{sys:ODE-Homo} 
   \dot{x} = F(x).
\end{equation}
The system is called (strictly) positive if $x \mapsto \phi(t,x)$ is (strictly) order-preserving for all $t>0$.
\end{definition}

Since we will work with continuous time systems. Theorem \ref{thm:main-pos} needs to be adapted for
system \eqref{sys:ODE-Homo}.

\begin{corollary} \label{coro:main}
Let $ \coneK$ be a normal cone in a finite dimensional vector
space $\vsV$.
Consider the dynamics \eqref{sys:ODE-Homo}. The following hold.
\begin{enumerate}
    \item[(i)] 
    If there exists some constant $a$, such that $aI + F $ is a strongly order-preserving mapping, 
    where $I$ is the identity mapping, then the system \eqref{sys:ODE-Homo} is strictly positive.
    \item[(ii)]
    If the system \eqref{sys:ODE-Homo} is strictly positive,
    then there exists a constant $\mu \in \mathbb{R}$, and $x_* \in \coneK \backslash \{0 \}$ (unique up
    positive scaling) such that  
    \begin{enumerate}
        \item $F(x_*) = \mu x_* $.
        \item $ d_{\mathcal{H}}(\phi(t,x), x_*) \to 0$ as $t \to \infty $, $\forall x\in \coneK \backslash \{0\}$.
        \item  If in addition, $F$ is linear, then the convergence is exponential, i.e., there exist
    two positive constants $k$, $\lambda$, such that
    $
        d_{\mathcal{H}} (\phi(t,x), x_*) \le k e^{-\lambda t } d_{\mathcal{H}}(x,x_*), \;
    $for all $t\ge 0$ and $x\in \coneK \backslash \{0\}$.
    \end{enumerate}

\end{enumerate}

\end{corollary}

\begin{pf}
(i) Rewrite the system \eqref{sys:ODE-Homo} as $\dot{x} = -ax + (aI+F)(x)$. 
Without loss of generality, assume $\mathcal{V}= \mathbb{R}^n$.
Then 
\begin{align*}
    \phi(t,x) = e^{-at} x + \int_0^t e^{-a(t-s)}(aI+F)(\phi(s,x)) ds
\end{align*}from which we deduce that $x\mapsto \phi(t,x)$ is strongly order-preserving for every $t>0$.

(ii) Since $x\mapsto \phi(t,x)$ is strongly order-preserving for all $t>0$,
by Theorem \ref{thm:main-pos}, 
for every fixed $t>0$, there exists some $x_*$ such that 
\begin{equation} \label{eq:phi: lambda}
\phi(t,x_*) = \lambda_t x_* 
\end{equation}
and that
$\phi(kt, x)$ converges to $x_*$ in Hilbert metric as $k\to \infty$ for all $x\in \coneK \backslash \{0\}$.
It is obvious that $x_*$ is independent of $t$, i.e., the equation holds for all $t>0$ for some $t$-dependent
parameters $\lambda_t$ and that $\phi(t,x)\to x_*$ in Hilbert metric as $t\to \infty$ for all $x\in \coneK \backslash \{0\}$. Now \eqref{eq:phi: lambda} simply says that the ray $\{ r x_*: r \geq 0  \}$ is invariant.
Thus the time derivative of $\phi(t,x_*)$ at $t=0$ results in $F(x_*)= \dot{\lambda}_0 x_*$. 
Hence we can take $\mu = \dot{\lambda}_0$ in (ii-a).

\end{pf}

\section{Lyapunov equations}\label{sec:lyap}

\subsection{Continuous-time Lyapunov equations}

The first equation that we are going to study is the so-called continuous-time Lyapunov equation:
\begin{equation}
A^{\top}P+PA+C^{\top}C=0\label{eq: CT-LYP}
\end{equation}
where $A\in\mathbb{R}^{n\times n}$, $C\in\mathbb{R}^{m\times n}$
are known matrices, and the matrix $P$ is the unknown variable to
be solved. It is well-known that:
\begin{fct} \label{fact:Lin_LYP-CT}
If $(A,C)$ is observable and $A$ is stable, there exists a unique
positive definite solution $P_{\#}$ to (\ref{eq: CT-LYP}). 
\end{fct}

Below, we adopt a positive systems approach to prove this fact. To illustrate how an algebraic equation relates to a positive dynamical system, consider a simple algebraic equation:
\begin{equation} \label{eq:lin-pos}
Ap + c=0
\end{equation}where $A$ is a Metzler matrix (all off-diagonal elements are non-negative), and $c$ is a non-negative vector. It is well-known that this 
equation has a unique strictly positive solution $p$ when $A$ is stable.

\begin{figure}[ht]
    \centering
    \includegraphics[scale=0.6]{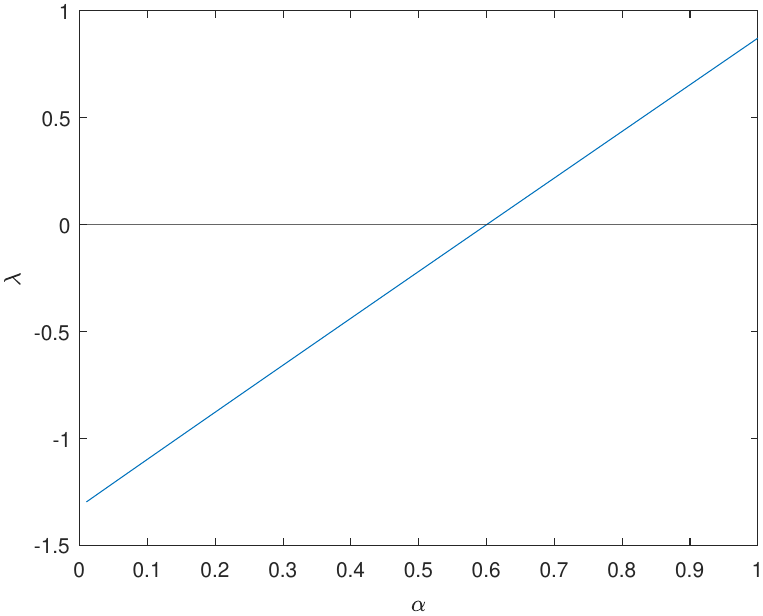}
    \caption{Eigenvalue corresponding to the Perron-Frobenius vector as a function of $\alpha$ for the matrix $A+\alpha c \allOne^\top$.}
    \label{fig:lin-pos-eq}
\end{figure}

To prove this by positive system theory, we may associate this equation with a dynamical system $\dot{x} = Ax+c$,
which is a positive system since $A$ is Metzler and $c\ge 0$ \cite{rantzer2018tutorial}.
However, the system is not linear (neither homogeneous). To apply Corollary \ref{coro:main}, we ``linearize'' the system in the following fashion:
\begin{equation*}
\dot{x}  =(A+\alpha c \allOne^{\top})x
\end{equation*}
where $\alpha>0$ is a positive constant, and $\allOne$ is a column vector with all ones.
This is a linear strictly positive system. In view of Corollary \ref{coro:main}, it converges
(exponentially in Hilbert metric) to the unique Perron-Frobenius vector $\xi>0$ and that $(A+\alpha c\allOne^{\top})\xi=\lambda\xi$ holds
for some $\lambda\in\mathbb{R}$.
Suppose that $A$ is Hurwitz, 
then for $\alpha$ sufficiently small,
$\lambda$ must be negative; for $\alpha$ large,
$\lambda$ is clearly positive. By continuity (of the largest eigenvalue), there exists some $\alpha>0$,
such that $\lambda=0$, i.e., $A\xi+\alpha c\allOne^{\top}\xi=0$, see Fig. \ref{fig:lin-pos-eq}. Now let $p=\frac{\xi}{\alpha(\allOne^{\top}\xi)}$, we see that 
$p$ is the desired strictly positive solution to (\ref{eq:lin-pos}). 

Now going back to the Lyapunov equation \eqref{eq: CT-LYP},
we mimic the above procedure to define the following linear dynamics on $\Sn$:
\begin{equation}
\dot{P}=A^{\top}P+PA+\alpha{\rm tr}(P)C^{\top}C\label{eq:CT-ODE-trace}
\end{equation}
where $\alpha$ is a positive constant.

\begin{lemma}
\label{lem:pos-ct}If $(A,C)$ is observable, then the system (\ref{eq:CT-ODE-trace})
is strictly positive on $\Sn$.
\end{lemma}

\begin{pf}
Given initial state $P_{0}\in\Sn\backslash\{0\}$, the solution \eqref{eq:CT-ODE-trace} satisfies
\[
P(t)=e^{A^{\top}t}P_{0}e^{At}+\alpha\int_{0}^{t}{\rm tr}(P(s))e^{A^{\top}(t-s)}C^{\top}Ce^{A(t-s)}{\rm d}s
\]
If there exists a non-zero vector $v\in \mathbb{R}^n$ such that $P(t)v=0$ for
some $t>0$, then $v^{\top}P(t)v=0$ and 
\[
{\rm tr}(P(s))Ce^{A(t-s)}v=0,\;\forall s\in[0,t].
\]
Since $P(s)\succeq e^{A^{\top}s}P_{0}e^{A^{\top}s}$, ${\rm tr}(P(s))$
must be positive since $P_{0}\ne0$. Thus $Ce^{As}v=0$ for all $s\in[0,t]$,
which is impossible thanks to the observability of $(A,C)$. Thus
$P(t) \succ 0$, or $P(t)\in{\rm Int}\,\Sn$ for all $t>0$. 
\end{pf}

With Lemma \ref{lem:pos-ct} at hand, we can now prove Fact \ref{fact:Lin_LYP-CT} using positive systems theory.

\textbf{Positive system argument for {Fact \ref{fact:Lin_LYP-CT}}: }
Since (\ref{eq:CT-ODE-trace}) is a strictly positive system
on $\Sn$, by Corollary \ref{coro:main}, there exists a unique
direction $P_{*} \succ 0$ and some constant $\lambda\in\mathbb{R}$, such
that
\[
A^{\top}P_{*}+P_{*}A+\alpha{\rm tr}(P_{*})C^{\top}C=\lambda P_{*}
\]
Let $v\in\mathbb{C}^{n}\backslash\{0\}$ be an eigenvector of $A$
and $\eta\in\mathbb{C}$ the corresponding eigenvalue, then
$\lambda=2{\rm Re}(\eta)+\alpha\frac{{\rm tr}(P_{*})||Cv||^{2}}{v^{*}P_{*}v}$.
Due to observability of $(A,C)$, $Cv\ne 0$, which yields $||Cv||^2 \ne 0$.
Since $A$ is stable,
we must have $\lambda<0$ when $\alpha$ is sufficiently small,
and when $\alpha $ becomes large, $\lambda$ must exceed zero.
By continuity of the largest eigenvalue,
there exists some $\alpha>0$ such that $\lambda=0$, i.e.,
\[
A^{\top}P_{*}+P_{*}A+\alpha{\rm tr}(P_{*})C^{\top}C=0
\]
but then $P_{\#}=\frac{P_{*}}{\alpha{\rm tr}(P_{*})}$ solves \eqref{eq: CT-LYP}.
Since $P_{*}$
is unique up to positive scaling,
$P_{\#}$ is unique. Furthermore, with this $\alpha $ fixed, the dynamics \eqref{eq:CT-ODE-trace} satisfy
\[
 d_{\mathcal{H}}(P(t),P_\#) \le ke^{-\lambda t} d_{\mathcal{H}}(P_0, P_\#), \, \forall t\ge 0, \, P_0 
 \succeq 0
\]for any some positive constants $k, \, \lambda$.

Similar arguments hold for discrete-time Lyapunov equations, as shown in the next subsection.

\subsection{Discrete-time Lyapunov equation}

In this subsection, we study the discrete counterpart of (\ref{eq: CT-LYP}),
i.e., the equation
\begin{equation}\label{eq:DT-LYP}
A^{\top}QA-Q+C^{\top}C=0.
\end{equation}

To analyze this equation from a fixed-point perspective, a natural approach is to consider the map $h(X)=A^\top XA + C^\top C$, as the solution to \eqref{eq:DT-LYP} corresponds to a fixed-point of $h$. However, applying fixed-point theorems typically requires certain structures, such as convexity and compactness. In \cite{pates2024lyapunov}, Pates considered an alternative mapping
\begin{equation}
f (X) = \frac{A^\top X A + \alpha C^\top C}{{\rm tr}(A^\top XA + \alpha C^\top C)},
\end{equation}which is invariant on the compact convex set $ \{ X\succeq 0: {\rm tr}(X) = 1 \}$. 
This invariance enables the application of Brouwer’s fixed-point theorem to guarantee the existence of a fixed-point. However, a known limitation of Brouwer’s theorem is that it does not provide any assurance of uniqueness.

It is worth noting that $f$ defined in this way is non-homogeneous and
thus positive systems theory is not directly applicable.
In what follows, we associate the equation \eqref{eq:DT-LYP} with a linear positive system, allowing us to apply Theorem \ref{thm:main-pos}. As a result, both existence and uniqueness of the solution are ensured. The goal is to prove the following fact:
\begin{fct} \label{fact:DT-LYP}
If $A$ is Schur stable and $(A,C)$ observable. Then there exists
a unique $Q  \succ 0$ solving the Lyapunov equation (\ref{eq:DT-LYP}).
\end{fct}

We consider the discrete-time dynamics:
\begin{equation} \label{sys:DT-LYP}
X_{k+1}=F(X_{k})
\end{equation}
on $\Sn$, where 
\[
F(X)=A^{\top}XA+\alpha{\rm tr}(X)C^{\top}C
\]
with $\alpha$ a positive constant. 

The following lemma shows that observability ensures {\it ultimate strict positivity} of the system \eqref{sys:DT-LYP} in the sense that there exists some $T>0$ such that
$X_k \succ 0$ for all $k\ge T$ and $X_1 \in {\rm Int} \; \Sn $.

\begin{lemma}\label{lem:ultimate-DT}
If $A$ is Schur stable and $(A,C)$ is observable, then the system \eqref{sys:DT-LYP} is ultimately strictly positive.
\end{lemma}
\begin{pf}
Let $X_{0} \in  {\Sn} \backslash \{0\}$, then $X_{k}=(A^{k})^{\top}X_{0}A^{k}+\sum_{j=0}^{k-1}\alpha{\rm tr}(X_{j})(A^{j})^{\top}C^{\top}CA^{j}$
for all $k\ge1$. If there exists $v$ such that $v^{\top}X_{k}v=0$,
then $CA^{j}v=0$ for all $j\ge0$, which implies $v=0$ thanks to the observability of $(A,C$). Hence $X_{k}\in{\rm Int}\,\Sn$ for all $k \ge n$
as claimed. 
\end{pf}

\begin{figure}[ht]
    \centering
    \includegraphics[scale=0.6]{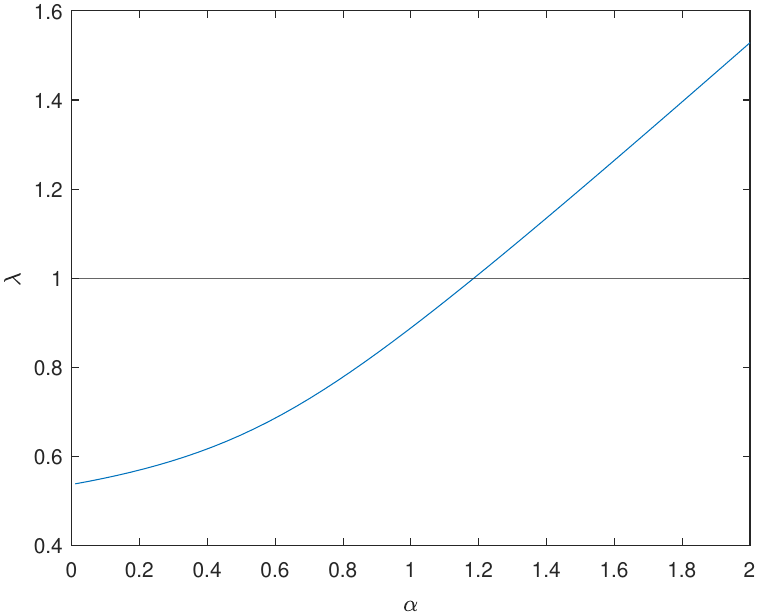}
    \caption{Eigenvalue corresponding to the Perron-Frobenius vector as a function of $\alpha$ for the operator $L(X) = A^\top X A+ \alpha {\rm tr}(X) C^\top C$ with $A$ Schur stable and $(A,C)$ observable.}
    \label{fig:DT-LYP-eig}
\end{figure}

\textbf{Positive system argument of Fact \ref{fact:DT-LYP}}:
In view of Lemma \ref{lem:ultimate-DT} and Theorem \ref{thm:main-pos}, there exists $\lambda>0$
and a unique (up to scaling) $X_{*}\in{\rm Int}\;\Sn$,
such that $F(X_{*})=\lambda X_{*}$. Note
\[
\lambda\ge\frac{\alpha{\rm tr}(C^{\top}C)}{{\rm tr}(X_{*})}\to\infty
\]
as $\alpha\to\infty$; and when $\alpha$ is sufficiently small, it is
clear that $\lambda<1$ since $A$ is Schur stable. By continuity,
there exists some $\alpha>0$, such that $\lambda=1$, i.e.,
\[
A^{\top}X_{*}A+\alpha{\rm tr}(X_{*})C^{\top}C=X_{*},
\]as shown in Fig. \ref{fig:DT-LYP-eig}.
Now $ Q :=\frac{X_{*}}{\alpha{\rm tr}(X_{*})}$ solves (\ref{eq:DT-LYP}) and that
there exists some positive constant $k$, $\mu$ such that
\[
    \dH (X_k ,Q) \le \mu e^{-\lambda k} \dH (X_1, Q)
\]for all $k\ge 1$ and $X\in \Sn \excdZero$.
Uniqueness is obvious.

\subsection{Discussion}

We summarize the major steps in deriving the solutions of Lyapunov equations. 
\begin{enumerate}
    \item[a)] Construct a linearized positive dynamical system with a tuning parameter $\alpha$ that controls the spectral radius of the positive operator.
    
    \item[b)] Show that the system is strictly positive based on additional assumptions, typically observability assumptions.
    
    \item[c)] By tuning the parameter $\alpha$, and making use of the continuity of the spectral radius as a function of $\alpha$, show that the spectral radius can reach a desired value.

    \item[d)] Show that the dynamical system gives a solution to the equation when it is at steady state.
\end{enumerate}

By following these steps, we can easily extend the preceding results to more general types of Lyapunov equations, such as 
\begin{equation}\label{eq:gen_LYP_stochastic}
A^\top P + PA + \sum_{i=1}^r A_i^\top P A_i + C^\top C=0,
\end{equation}
and 
\begin{equation}
    \label{eq:gen_LYP}
    A^\top P E + E^\top PA + C^\top C = 0.
\end{equation}
Equation \eqref{eq:gen_LYP_stochastic} arises naturally in linear stochastic systems \cite{brockett1973lie} and equation \eqref{eq:gen_LYP} appears in singular linear systems \cite{syrmos1995discrete}, see also \cite{benner2011lyapunov} and the references therein. 
For example, by defining systems
\begin{equation*}
\begin{aligned}
\dot{X} & =A^{\top}X+XA+\sum_{j=1}^{m}A_{j}^{\top}XA_{j}+\alpha{\rm tr}(X)C^{\top}C \\
\dot{X} & = A^\top X E + E^\top XA + \alpha {\rm tr}(X)C^\top C
\end{aligned}
\end{equation*}
one can easily establish the solution properties of the equations  \eqref{eq:gen_LYP_stochastic} and \eqref{eq:gen_LYP}.

\section{Algebraic Riccati equations}\label{sec:ricc}

In this section, we explore to what extent we can extend the previous results to Riccati equations. The main difficulty, as shall be seen, arises from the fact that Riccati equations cannot be ``linearized'' in the same way as Lyapunov equations. Nevertheless, Riccati equations can be ``homogenized'', allowing us to apply Theorem \ref{thm:main-pos} and Corollary \ref{coro:main}.

\subsection{Discrete-time Algebraic Riccati Equation} 

We start with the discrete-time algebraic Riccati equation:
\begin{equation}
P=A^{\top}PA-A^{\top}PB(R+B^{\top}PB)^{-1}B^{\top}PA+C^{\top}C\label{eq:Riccati}
\end{equation}
where $R \succ 0$. 
The goal in this subsection is to provide a positive systems perspective for the following classical result:
\begin{fct} \label{fact:Ric-DT}
If $(A,B)$ is controllable and $(A,C)$ is observable, then there exists
a unique positive definite solution to (\ref{eq:Riccati}). 
\end{fct}

Define the system
\begin{equation}
X_{k+1}=F(X_{k})\label{sys:DARE}
\end{equation}
with
\begin{equation} \label{eq:F-DT_Riccati}
F(X) = A^\top S(X) A + \alpha {\rm tr }(X) C^\top C
\end{equation}for a positive constant $\alpha$, where 
\[ 
S(X):=X-XB(\alpha{\rm tr}(X)R+B^{\top}XB)^{-1}B^{\top}X.
\] 

\begin{remark}
The introduction
of the scalar-valued function $ P\mapsto {\rm tr}(P) $ is instrumental 
for the application of positive system theory,
as it renders the system \eqref{sys:DARE} homogeneous.
\end{remark}

\begin{remark}
    The system \eqref{sys:DARE} resembles the well-known value iteration
    \[
        P_{k+1} = A^{\top}P_k A-A^{\top}P_k B(R+B^{\top}P_k B)^{-1}B^{\top}P_k A+C^{\top}C,
    \]see, for example \cite{bertsekas2012dynamic,lee2022tac,lai2025robust}.
    The key difference between \eqref{sys:DARE} and value iteration is that
    value iteration is not homogeneous. As a consequence, their convergence properties
    differ. 
    More importantly, by ``homogenizing'' the value iteration procedure, positive
    systems theory can immediately be applied to draw conclusions on the convergence of the iteration and the existence and uniqueness of the corresponding solution.
    
\end{remark}

\begin{lemma} \label{lem:Ric-DT}
The operator $F$ defined in \eqref{eq:F-DT_Riccati} is 
continuous, homogeneous, order-preserving, and concave. 
If $(A,C)$ is observable, then $F^n$ is strongly order-preserving.
\end{lemma}
\begin{pf}
Continuity and homogeneity of $F$ are obvious; order-preservation follows from
Lemma \ref{lem:S(X)} (b); the concavity of the Schur complement 
is well-known, for a reference, see \cite{ando1979concavity}. 

Assume now $(A,C)$ is observable. Let $X_1 \in \Sn \backslash \{ 0 \}$.
Given $k>0$, if there exists some $v\ne0$ such that $v^{\top}X_{k+1}v=0$,
then 
\[
v^{\top}X_{k+1}v=v^{\top}(A^{\top}S(X_{k})A+\alpha{\rm tr}(X_{k})C^{\top}C)v=0
\]
which implies $Cv=0$, and $v^{\top}A^{\top}S(X_{k})Av=0$. In view
of Lemma \ref{lem:S(X)} (c), the latter is equivalent to $v^{\top}A^{\top}X_{k}Av=0$.
But 
\[
\begin{aligned}
v^{\top}A^{\top}X_{k}Av   =v^{\top}A^{\top}[(A^{\top}S(X_{k-1})A \\ 
  \quad +\alpha{\rm tr}(X_{k-1})C^{\top}C)]Av=0.
\end{aligned}
\]
implying $CAv=0$ and $v^{\top}(A^{2})^{\top}X_{k-1}A^{2}v=0$. Proceeding
this recursively, we arrive at 
\[
{\rm col}(C, CA,\cdots , CA^{k-1}) v = 0
\]
Now choosing $k=n$ yields $v=0$. Hence we arrive at the following implication
\begin{equation} \label{eq:X to Fn}
X\in \Sn \backslash \{0 \} \Rightarrow F^n (X) \in {\rm Int} \; \Sn
\end{equation}
This implies $F^n$ is strongly order-preserving. Indeed, let $Y \succeq X \succeq 0$ and $Y \neq X$,
then
\begin{align*}
    F^n (Y) & = F^n( Y-X + X) \\
    & \succeq F^n(Y-X) + F^n(X)  \\
    & \succ F^n(X).
\end{align*}where we have used the fact that $F^n(Y-X) \succ 0$ and 
that concavity and homogeneity together imply
super-additivity\footnote{A mapping $F$ on cone $\coneK$ is called super-additive if $F(x+y) \ge_{\coneK} F(x) + F(y)$.}.
\end{pf}

\begin{remark}
    Note that for any $m \in \mathbb{N}$, the power $F^m$ is also continuous, order-preserving, and concave (composition of two monotone increasing concave maps is again concave). 
\end{remark}

We are now ready to prove Fact \ref{fact:Ric-DT} from a positive systems perspective.

\textbf{Positive system argument for Fact \ref{fact:Ric-DT}}:
By Lemma \ref{lem:Ric-DT} and Theorem \ref{thm:main-pos},
there exists an $X_* \succeq 0 $ such that $F(X_{*})=\rho(F)X_{*}$. But then
$F^n (X_*) = \rho(F)^n X_*$. Since $F^n$ is strongly-order preserving, we must have 
$X_* \succ 0$, which is unique up to positive scaling.
Rewrite $F(X_*) = \rho(F) X_*$:
\[
\begin{aligned}
A^{\top}X_{*}A  -A^{\top}X_{*}B(\alpha  {\rm tr}  ( & X_{*}  )R   +B^{\top}X_{*}B)^{-1}B^{\top}X_{*}A \\ 
& +\alpha{\rm tr}(X_{*})C^{\top}C=\rho(F)X_{*}.
\end{aligned}
\]
Note that $F(X_{*})\preceq  A^{\top}X_{*}A+\alpha{\rm tr}(X_{*})C^{\top}C$;
hence if $A$ is Schur stable, there holds $\rho(F)<1$ when $\alpha$
is sufficiently small. 
By assumption, $(A,B)$ is controllable, then
there exists $K$ such that $A_{c}=A-BK$ is Schur stable. It is readily
checked that
\[
\begin{aligned}
& A^{\top}[X-XB(B^{\top}XB)^{-1}  B^{\top}X]A \\
 =& A_{c}^{\top}[X-XB(B^{\top}XB)^{-1}B^{\top}X]A_{c}
\end{aligned}
\]
Thus without loss of generality, we may assume that $A$ is Schur
stable, which implies that $\lambda<1$ for sufficiently small $\alpha>0$.

On the other hand, $\rho(F)\to\infty$ as $\alpha\to\infty$
due to observability of $(A,C)$. By continuity of the spectral
radius \cite[Corollary 5.5.5]{lemmens2012nonlinear}, there exists some
$\alpha>0$, and $X_{*} \succ 0$ such that $F(X_{*})=X_{*}$, i.e., 
\begin{align*}
A^{\top}X_{*}A  -A^{\top}X_{*}B(\alpha{\rm tr}(X_{*})R 
+B^{\top}X_{*}B)^{-1}B^{\top}X_{*}A \\
 +\alpha{\rm tr}(X_{*})C^{\top}C=X_{*}.
\end{align*}
Then $P =\frac{X_{*}}{\alpha{\rm tr}(X_{*})}$ solves (\ref{eq:Riccati}), which is obviously
the unique solution. Moreover, 
\[
d_{\mathcal{H}} (X_k, P) \to 0, \text{ as } k\to \infty, \forall X_1 \in \Sn \backslash \{0 \}.
\]

\subsection{Continuous-time Algebraic Riccati Equation}

Finally, study the continuous-time algebraic Riccati equation:
\begin{equation} \label{eq:CARE}
A^{\top}P+PA-PBR^{-1}B^{\top}P+C^{\top}C=0.
\end{equation}

The goal is to prove the following fact by positive systems theory:
\begin{fct} \label{fact:Ric-CT}
If $(A,B)$ is controllable and $(A,C)$ is observable, then there exists
a unique positive definite solution to (\ref{eq:CARE}). 
\end{fct}

To study this equation, we propose the following homogeneous system:
\begin{equation*}
\dot{P}=\alpha A^{\top}P+\alpha PA-\frac{1}{{\rm tr}(P)}PBR^{-1}B^{\top}P+\alpha^{2}{\rm tr}(P)C^{\top}C
\end{equation*}
with $P(0)=P_{0}$ and $\alpha $ a positive constant. Equivalently, the system can be written as
\begin{equation} \label{eq:sys:Riccati-ODE}
\dot{P}=A_{\alpha}(t) ^\top P +P A_{\alpha}(t)+\alpha^{2}{\rm tr}(P)C^{\top}C
\end{equation}
where $A_{\alpha}(t)=\alpha A-\frac{1}{2{\rm tr}(P)}BR^{-1}B^{\top}P(t)$.

\begin{lemma} \label{lem:Ric-CT}
If $(A,C)$ is observable, then the system \eqref{eq:sys:Riccati-ODE} is strictly positive.
\end{lemma}
\begin{pf}
Step 1: Show that
\[
P(0)=P_0 \in \Sn \excdZero \Rightarrow P(t) \in {\rm Int} \; \Sn, \; \forall t>0.
\]
Let $P(0)=P_0 \in \Sn \backslash \{0\}$ and $\Phi(t,\tau)$ be the solution to $\dot{X}=A_{\alpha}(t)X$, $X(\tau)=I$.
Then
\begin{align*}
P(t) & =\Phi(t,\tau)^{\top}P(\tau)\Phi(t,\tau) \\ 
& \quad +\alpha^{2}\int_{\tau}^{t}{\rm tr}(P(s))\Phi(t,s)^{\top}C^{\top}C\Phi(t,s){\rm d}s
\end{align*}
for all $\tau\in[0,t]$. If there exists some $t>0$ such that $P(t)$
is not positive definite, then there exists $v\ne0$, such that 1)
$P(\tau)\Phi(t,\tau)v=0$; 2) $C\Phi(t,\tau)v=0$ for all $\tau\in[0,t]$.
Differentiating 2) with respect to $\tau$ yields $CA_{\alpha}(\tau)\Phi(t,\tau)v=0$,
from which we deduce $CA\Phi(t,\tau)v=0$ for all $\tau\in[0,t]$
thanks to 1). Repeating this procedure and setting $\tau$ to $t$,
we arrive at
\[
{\rm col}(C, CA,\cdots , CA^{n-1}) v = 0
\]
Thus $v=0$ since $(A,C)$ is observable, from which it follows that $P(t) \succ 0$ for all $t>0$. 

Step 2: Show that $\phi(t, X+Y) \ge \phi(t,X) + \phi(t,Y)$ for all
$X, Y \in \Sn$. This can be done by proving the same inequality for 
the Euler approximation 
$
T_\tau (X) : = X + \tau(\alpha A^\top X + \alpha XA - XBR^{-1}B^\top X / {\rm tr}(X) +
\alpha^2 {\rm tr}(X)C^\top C $
when $\tau>0 $ is sufficiently small and then take the limit. But $T_\tau $ is concave and 
homogeneous, therefore it is also super-additive.

Finally, Step 1 and Step 2 together imply that 
\[
Y \succeq X \text{ and } Y \ne X \Rightarrow \phi(t,Y) \succ \phi(t,X)
\]for all $t>0$. The proof is complete.

\end{pf}

\textbf{Positive system argument for Fact \ref{fact:Ric-CT}:}  Since the system \eqref{eq:sys:Riccati-ODE}
is strictly positive by Lemma \ref{lem:Ric-CT}, by Corollary \ref{coro:main}, there exists a constant $\mu$,
and a unique (up to positive scaling) $X_* \succ 0$ such that
\begin{equation} \label{eq:CT-Ric-alpha}
\begin{aligned}
\alpha A^{\top}X_{*}+\alpha X_{*}A-\frac{1}{{\rm tr}(X_{*})}X_{*}BR^{-1}B^{\top}X_{*}  \\ +\alpha^{2}{\rm tr}(X_{*})C^{\top}C  
 = \mu X_{*}
\end{aligned}
\end{equation}
First, we show that $\mu<0$ when $\alpha$
is sufficiently small. For that, rewrite \eqref{eq:CT-Ric-alpha}
as
\begin{equation} \label{eq:dual-CT-Ric-alpha}
\begin{aligned}
\alpha Y_* A^{\top}+\alpha A Y_* -\frac{1}{{\rm tr}(X_{*})}BR^{-1}B^{\top}  \\
+\alpha^{2}{\rm tr}(X_{*}) Y_* C^{\top}C Y_*   
 = \mu Y_{*}
\end{aligned}
\end{equation}where $Y_*  = X_*^{-1}$. Let $v \in \mathbb{C}^n \excdZero$ be a 
left eigenvector of $A$, i.e., $v^* A = \eta v^*$ for some $\eta \in \mathbb{C}$. 
Then 
\begin{equation*} 
\begin{aligned}
2\alpha {\rm Re}(\eta) v^*Y_* v -\frac{1}{{\rm tr}(X_{*})}v^*BR^{-1}B^{\top}v  \\
+\alpha^{2}{\rm tr}(X_{*}) v^* Y_* C^{\top}C Y_* v  
 = \mu v^* Y_{*} v.
\end{aligned}
\end{equation*}Now that $(A,B)$ is controllable, we have $B^\top v \ne 0$, and as a 
result, the term $v^* BR^{-1} B^\top v$ is strictly positive which dominates the
terms on the left-hand side of the equation. Hence $\mu <0$ when $\alpha $ is small.

Similarly, one can show that $\mu>0$ when $\alpha $ is large due to the observability
of $(A,C)$.
Now by the continuity of cone spectral radius, there exists
some $\alpha>0$ such that $\mu=0$. This concludes our proof by choosing
$P = \frac{X_*}{\alpha^2 {\rm tr}(X_*)}$, the unique positive definite solution to \eqref{eq:CARE}. And 
\[
\dH (P(t), P) \to 0 \text{ as } t\to \infty, \forall P(0)\in \Sn \excdZero.
\]

\begin{remark}

We observe that in establishing the spectral properties of the dynamical systems associated with algebraic Riccati equations -- both in discrete and continuous time -- controllability and observability play equally essential roles. Controllability contributes to reducing the spectral radius, while observability tends to increase it. Crucially, the absence of either property prevents us from guaranteeing the existence or uniqueness of solutions.

\end{remark}

\subsection{Riccati equations in positive control systems}
In \cite{rantzer2022explicit}, Rantzer considered the optimal control problem of a linear positive system with linear cost, see also
\cite{blanchini2023optimal}
for the continuous-time counterpart. The Bellman equation -- or algebraic Riccati equation -- associated with this problem was derived:
\begin{equation} \label{eq:Ric-positive}
    p = s + A^\top p - E^\top |r + B^\top p|
\end{equation}where $s$, $A$, $E$, $r$ and $B$ are known matrices or vectors.
It is assumed that $E \ge 0$ (point-wise non-negative), $s \ge E^\top r$, and $A \ge |BE| $. Here $|\cdot|$ stands for point-wise absolute value.

Li and Rantzer studied the solution properties of the equation \eqref{eq:Ric-positive} in \cite{li2024exact} under certain assumptions, similar to the controllability
and observability assumptions made in LQR control:
\begin{enumerate}
    \item[{\bf H1}] there exists a matrix $L$, with $|L|\le E$, such that $A-BL$ is Schur stable;
    \item[{\bf H2}] there holds:
    \[
    (s-E^\top |r|)^\top \sum_{i=0}^{n-1} (A-|B|E)^i >0
    \]
\end{enumerate}

As before, we adopt a positive systems approach to study the equation \eqref{eq:Ric-positive}.
Consider the discrete-time linear homogeneous system
\begin{equation}
    p_{k+1} = F(p_k)
\end{equation} where
\[
F(p) = \alpha (\mathds{1}^\top p) s + A^\top p - E^\top |\alpha(\mathds{1}^\top p) r 
+ B^\top p|
\]for some positive constant $\alpha$.
Notice that $F$ is concave, homogeneous, and that $F^n$ satisfies (by {\bf H2})
\[
p \in \mathbb{R}^n_+ \excdZero \Rightarrow F^n (p) \in {\rm Int } \; \mathbb{R}^n_+.
\] Thus $F$ is sup-additive, implying that $F^n$ is strictly order preserving. As a consequence, Theorem \ref{thm:main-pos} implies the existence of a unique (up to positive scaling) $p_*\in {\rm Int} \; \mathbb{R}^n_+$, such that 
$F(p_*) = \lambda p_*$ for some positive constant $\lambda$. Next, we shall show that $\lambda$ 
can be taken as $1$ by adjusting the parameter $\alpha$. Indeed, when $\alpha =0$,
\begin{align*}
    \lambda p_*  & = A^\top p_* - E^\top |B^\top p_*|  \\
    & \le A^\top p_* - L^\top B^\top p_* \\
    & = (A-BL)^\top p_*,
\end{align*}where we have used the assumption $|L|\le E$ in the second line. But $A-BL$ is positive and Schur stable by {\bf H1}, we must have $\lambda <1$.
By continuity, $\lambda <1$ for all sufficiently small $\alpha>0$. On the other hand, 
by {\bf H2}, we have
\begin{align*}
    \lambda^n p_* & = F^n(p_*) \\
    & > \alpha(\mathds{1}^\top p_*) (s-E^\top |r|)^\top \sum_{i=0}^{n-1} (A-|B|E)^i p_* \\
    & > 0.
\end{align*} Hence $\lambda \to \infty$ when $\alpha \to \infty$, as desired.
The existence and uniqueness of the solution to \eqref{eq:Ric-positive} are now justified.

\section{Conclusion} \label{sec:conclude}
We have provided new insights into both Lyapunov and algebraic Riccati equations through the lens of positive systems theory. These findings establish interesting connections between algebraic equations and dynamical systems, demonstrating the efficacy of positive system theory. 
Future research may explore time-varying equations and quantify error bounds of the iteration methods proposed in the paper.

\section{Acknowledgement}
We have received constructive suggestions from Shiyong Zhu during the preparation
of the manuscript.

\section*{Appendix}

\begin{lemma}[\cite{anderson1971shorted}]
\label{lem:mono-Schur}
    The Schur complement is monotone in the following sense.
    Let $P, Q$ be two positive semi-definite matrices 
    \[
        P = 
        \begin{bmatrix}
            P_{1} & P_{2} \\ 
            P_{2}^\top & P_3
        \end{bmatrix},\quad
        Q = 
        \begin{bmatrix}
            Q_1 & Q_2 \\
            Q_2^\top & Q_3
        \end{bmatrix}.
    \]Then $P\succeq Q$ implies 
    \[
         P_1 - P_2 P_3^\dagger P_2^\top \succeq Q_1 - Q_2 Q_3^\dagger Q_2^\top \succeq 0.
    \]
\end{lemma}

\begin{lemma}[\cite{DeSouza1989tac}]
\label{lem:S(X)}Given matrices $X\succeq 0$, $R\succ 0$ and $B$ a matrix with appropiate dimension. Let $S(X)=X-XB(R+B^{\top}XB)^{-1}B^{\top}X$, then
\begin{enumerate}
    \item[(a)] $S(X)\succeq 0 $;
    \item[(b)] $X_1 \succeq X_2 \succeq 0$ implies $S(X_1) \succeq S(X_2) \succeq 0$.
    \item[(c)] $\ker X=\ker S(X)$;
\end{enumerate}
\end{lemma}

\begin{pf}
(a) and (b) are results from monotonicity of the Riccati operator,
see for example \cite[Lemma 3.1]{DeSouza1989tac}.
We provide a proof for (c). Let $v$ be such that $v^{\top}S(X)v=0$,
or $v^{\top}Xv-v^{\top}XB(R+B^{\top}XB)^{-1}B^{\top}Xv=0$. Note
\[
\begin{bmatrix}v^{\top}Xv & v^{\top}XB\\
B^{\top}Xv & R+B^{\top}XB
\end{bmatrix}\succeq 0
\]
If $v^{\top}Xv\ne0$, this matrix must be positive definite since
\[
R+B^{\top}XB-B^{\top}Xv(v^{\top}Xv)^{-1}v^{\top}XB\succeq R\succ 0,
\]
leading to a contradiction. Thus $Xv =0$, implying that $\ker S(X)\subseteq\ker X$. The
reverse inclusion is obvious. %

\end{pf}

\bibliographystyle{plain}
\bibliography{mono}

\end{document}